\documentclass[a4paper, twoside, 12pt]{article}
\usepackage{mathrsfs}
\usepackage{amsfonts}
\usepackage{amsmath}
\usepackage{amssymb}
\usepackage{amsthm}
\usepackage{hyperref}
\usepackage{cite}

 \numberwithin{equation}{section}
\theoremstyle{plain}\newtheorem{thm}{Theorem}[section]
\theoremstyle{plain}\newtheorem{lem}{Lemma}[section]
\theoremstyle{definition}
\theoremstyle{definition}
\theoremstyle{plain}
\theoremstyle{plain}\newtheorem{corr}{Corollary}[section]
\theoremstyle{plain}\newtheorem{rem}{Remark}[section]
\newtheorem*{con}{Conjecture of P$\acute{\textnormal{o}}$lya}
\newtheorem*{ack}{Acknowledgment}

\usepackage{fancyhdr}
\date{}

\title { Lower Bounds for Laplacian and\\ Fractional Laplacian Eigenvalues}
\author{Guoxin Wei \ \ He-Jun Sun\ \ {\rm and}\ \ Lingzhong Zeng}
\begin{document}
\maketitle

{\narrower\noindent \small { \textsc{Abstract:}}\ \ In this paper,
we investigate eigenvalues of Laplacian on a bounded domain in an
$n$-dimensional Euclidean space and obtain a sharper lower bound for
the sum of its eigenvalues, which gives an improvement of results
due to A. D. Melas \cite{M}. On the other hand, for the case of
fractional Laplacian $(-\Delta)^{\alpha/2}|_{D}$, where
$\alpha\in(0,2]$, we obtain a sharper lower bound for the sum of its
eigenvalues, which gives an improvement of results due to S.Y. Yolcu
and T. Yolcu \cite{YY}. }

\footnotetext{



The first author and the second author were supported by the
National Natural Science Foundation of China (Grant Nos.11001087,
11001130). }

\section{Introduction}

Let $D\subset \mathbb{R}^{n}$ be a bounded domain with piecewise
smooth boundary $\partial D$ in an $n$-dimensional Euclidean space
$\mathbb{R}^{n}$. Let $\lambda_{i}$ be the $i$-th eigenvalue of the
fixed membrane problem:
\begin{equation}\begin{cases}\label{1.1}\Delta u+\lambda u=0,~~~~~~~~~~~~~in~D,\\
u=0,~~~~~~~~~~~~~~~~~~~~~on~\partial D,\end{cases}\end{equation}
where $\Delta$ is the Laplacian in $\mathbb{R}^{n}$. It is well
known that the spectrum of this eigenvalue problem is real and
discrete:
$$0<\lambda_{1}\leq\lambda_{2}\leq\lambda_{3}\leq\cdots\rightarrow+\infty,$$ where each $\lambda_{i}$ has finite
multiplicity which is repeated according to its multiplicity. If we
use the notations $Vol(D)$ and $\omega_{n}$ to denote the volume of
$D$  and the volume of the unit ball in $\mathbb{R}^{n}$ ,
respectively, then Weyl's asymptotic formula asserts that the
eigenvalues of the fixed membrane problem (\ref{1.1}) satisfy the
following formula:
\begin{equation}\label{1.2}\lambda_{k}\sim \frac{4\pi^{2}}{(\omega_{n}Vol(D))^{\frac{2}{n}}}k^{\frac{2}{n}},~~k \rightarrow+\infty.\end{equation}
From the above asymptotic formula, it follows directly that
\begin{equation}\label{1.3}\frac{1}{k}\sum^{k}_{i=1}\lambda_{i}\sim\frac{n}{n+ 2}\frac{4\pi^{2}}
{(\omega_{n}Vol(D))^{\frac{2}{n}}}k^{\frac{2}{n}},~~k
\rightarrow+\infty.
\end{equation}\label{1.4} P$\acute{\textnormal{o}}$lya \cite{Pol} proved that \begin{equation}\lambda_{k}\geq \frac{4\pi^{2}}
{(\omega_{n}Vol(D))^{\frac{2}{n}}}k^{\frac{2}{n}},~~\textnormal{for}~k
=1,2,\cdots,\end{equation} if $D$ is a tiling domain in
$\mathbb{R}^{n}$ . Furthermore, he put forward the following:

\begin{con}If $D$ is a bounded domain in $\mathbb{R}^{n}$, then the $k$-th
eigenvalue $\lambda_{k}$ of the fixed membrane problem satisfies
\begin{equation}\label{1.5}\lambda_{k}\geq \frac{4\pi^{2}}
{(\omega_{n}Vol(D))^{\frac{2}{n}}}k^{\frac{2}{n}},~~\textnormal{for}~k
=1,2,\cdots.\end{equation}\end{con} On the Conjecture of
P$\acute{\textnormal{o}}$lya, Berezin \cite{Be} and Lieb \cite{Lie}
gave a partial solution.  In particular, Li and Yau [13] proved the
Berezin-Li-Yau inequality as follows:
\begin{equation}\label{1.6}\frac{1}{k}\sum^{k}_{i=1}\lambda_{i}\geq\frac{n}{n+
2}\frac{4\pi^{2}}
{(\omega_{n}Vol(D))^{\frac{2}{n}}}k^{\frac{2}{n}},~~\textnormal{for}~k
=1,2,\cdots.
\end{equation}
The formula (\ref{1.3}) shows that the result of Li and Yau is sharp
in the sense of average. From this inequality (\ref{1.6}), one can
derive
\begin{equation}\label{1.7}\lambda_{k}\geq\frac{n}{n+
2}\frac{4\pi^{2}}
{(\omega_{n}Vol(D))^{\frac{2}{n}}}k^{\frac{2}{n}},~~\textnormal{for}~k
=1,2,\cdots,
\end{equation} which gives a partial solution for the
conjecture of P$\acute{\textnormal{o}}$lya with a factor
$\dfrac{n}{n+2}$. We prefer to call this inequality (\ref{1.6}) as
Berezin-Li-Yau inequality instead of Li-Yau inequality because
(\ref{1.6}) can be obtained by a Legendre transform of an earlier
result by Berezin \cite{Be} as it is mentioned \cite{LW}. Recently,
improvements to the Berezin-Li-Yau inequality given by (\ref{1.6})
for the fixed membrane problem have appeared, for example see
\cite{KVW,M,W}. In particular, A.D.Melas \cite{M} has improved the
estimate (\ref{1.6}) to the following:
\begin{equation}\label{1.8}\frac{1}{k}\sum^{k}_{i=1}\lambda_{i}\geq\frac{n}{n+
2}\frac{4\pi^{2}}
{(\omega_{n}Vol(D))^{\frac{2}{n}}}k^{\frac{2}{n}}+\frac{1}{24(n+2)}\frac{Vol(D)}{Ine(D)},~~\textnormal{for}~k
=1,2,\cdots,
\end{equation}
where
$$Ine(D)=:\min_{a\in\mathbb{R}^{n}}\int_{D}|x-a|^{2}dx$$ is
called \emph{the moment of inertia} of $D$. After a translation of
the origin, we can assume that the center of mass is the origin
and$$ Ine (D) = \int_{D}|x|^{2}dx.$$

By taking a value nearby the extreme point of the function $f(\tau)$
(given by (\ref{3.5'})), we add one term of lower order of
$k^{-\frac{2}{n}}$ to its right hand side, which means that we
obtain a sharper result than (\ref{1.8}). In fact, we prove the
following:
\begin{thm}\label{thm1.1} Let $D$ be a bounded domain in an n-dimensional
Euclidean space $\mathbb{R}^{n}$. Assume that $\lambda_{i}, i =
1,2,\cdots,$ is the $i$-th eigenvalue of the eigenvalue problem
(1.1). Then the sum of its eigenvalues satisfies
\begin{eqnarray}\begin{aligned}\label{1.9}\frac{1}{k}\sum_{j=1}^{k}\lambda_{j}&\geq
\frac{nk^{\frac{2}{n}}}{n+2}\omega^{-\frac{2}{n}}_{n}(2\pi)^{2}Vol(D)^{-\frac{2}{n}}+\frac{1}{24(n+2)}\frac{Vol(D)}{Ine(D)}
\\&\quad~+\frac{nk^{-\frac{2}{n}}}{2304(n+2)^{2}}\omega_{n}^{\frac{2}{n}}(2\pi)^{-2}\Bigg{(}\frac{Vol(D)}{Ine(D)}\Bigg{)}^{2}Vol(D)^{\frac{2}{n}}.\end{aligned}\end{eqnarray}
\end{thm}

Furthermore, we consider the fractional Laplacian operators
restricted to $D$,  and denote them by $(-\Delta)^{\alpha/2}|_{D}$,
where $\alpha\in(0,2]$. This fractional Laplacian can be defined by
$$(-\Delta)^{\alpha/2}u(x)=:\textbf{P.V.}\int_{\mathbb{R}^{n}}\frac{u(x)-u(y)}{|x-y|^{n+\alpha}}dy,$$ where
$\textbf{P.V.}$ denotes the principal value and
$u:~\mathbb{R}^{n}\rightarrow\mathbb{R}.$ Define the characteristic
function $\chi_{D}:~t\mapsto\chi_{D}(t)$ by
\begin{equation*}\chi_{D}(t)=\begin{cases}1,~~~~~~~~~~~~x\in D,\\
0,~~~~~~~x\in\mathbb{R}^{n}\backslash D,
\end{cases}\end{equation*}
then the special pseudo-differential operator can be represented as
the Fourier transform of the function $u$ \cite{Lan,V}, namely
\begin{equation*}(-\Delta)^{\alpha/2}|_{D}u:=\mathscr{F}^{-1}[|\xi|^{\alpha}\mathscr{F}[u\chi_{D}]],\end{equation*} where $\mathscr{F}[u]$ denotes the Fourier transform of a
function $u: \mathbb{R}^{n}\rightarrow\mathbb{R}$:
$$\mathscr{F}[u](\xi)=\widehat{u}(\xi)=\frac{1}{(2\pi)^{n}}\int_{\mathbb{R}^{n}}e^{-ix\cdot\xi}u(x)dx.$$

It is well known that the fractional Laplacian operator
$(-\Delta)^{\alpha/2}$ can be considered as the infinitesimal
generator of the symmetric $\alpha$-stable process
\cite{BG,BK1,BK2,BKS,YY}. Suppose that a stochastic process
$\mathcal {X}_{t}$ has stationary independent increments and its
transition density (i.e., convolution kernel)
$p^{\alpha}(t,x,y)=p^{\alpha}(t,x-y),~t>0,~x,y \in \mathbb{R}^{n}$
is determined by the following Fourier transform
$${\rm Exp}(-t |\xi|^{\alpha}) = \int_{\mathbb{R}^{n}}e^{i\xi\cdot
y}p ^{\alpha}(t,y)dy,~~t>0,~~\xi\in \mathbb{R}^{n},$$then we can say
that the process $\mathcal {X}_{t}$ is an $n$-dimensional
\emph{symmetric $\alpha$-stable process} with order $\alpha\in(0,
2]$ in $\mathbb{R}^{n}$(also see \cite{BK1,BK2,YY}).
\begin{rem}Given $\alpha=1,~\mathcal {X}_{t}$  is the Cauchy process in $\mathbb{R}^{n}$
whose transition densities are given by the Cauchy distribution
(Poisson kernel)$$p^{1}(t,x,y)=\frac{c_{n}t}{(t^{2}+|x
-y|^{2})^{\frac{n+1}{2}}},~~t>0,~~x,y\in\mathbb{R}^{n} ,$$ where
$$c_{n}=\Gamma(\frac{n+1}{2})/\pi^{\frac{n+1}{2}}=\frac{1}{\sqrt{\pi}\omega_{n}},$$ is the semiclassical constant
that appears in the Weyl estimate for the eigenvalues of the
Laplacian.
\end{rem}

\begin{rem}Given $\alpha= 2$, $\mathcal {X}_{t}$ is just the usual
$n$-dimensional Brownian motion $\mathcal {B}_{t}$ but running at
twice the speed, which is equivalent to say that, when $\alpha=2$,
we have $\mathcal {X}_{t}=\mathcal {B}_{2t}$  and
$$p^{2}(t,x,y)=\frac{1}{(4\pi t)^{n/2}}{\rm Exp}\Bigg{[}\frac{-|x-y|^{2}}{4t}\Bigg{]},~~t
>0,~~x,y\in\mathbb{R}^{n}.$$\end{rem}Let $\Lambda_{j}^{\alpha}$ and $u_{j}^{\alpha}$ denote the $j$-th
eigenvalue and the corresponding normalized eigenvector of
$(-\Delta)^{\alpha/2}|_{\Omega}$, respectively. Eigenvalues
$\Lambda_{j}^{\alpha}$ (including multiplicities) satisfy
$$0<\Lambda_{1}^{(\alpha)}\leq\Lambda_{2}^{(\alpha)}\leq\Lambda_{3}^{(\alpha)}\leq\cdots\rightarrow+\infty.$$

For the case of $\alpha=1$, E. Harrell and S. Y. Yolcu gave an
analogue of the Berezin-Li-Yau type inequality for the eigenvalues
of the Klein-Gordon operators $\mathscr{H}_{0,D}:=\sqrt{-\Delta}$
restricted to $D$ in \cite{HY}:
\begin{eqnarray}\label{1.10}\begin{aligned}\frac{1}{k}\sum^{k}_{j=1}\Lambda^{(\alpha)}_{j}\geq\frac{n}{n+1}
\Bigg{(}\frac{2\pi}{(\omega_{n}Vol(D))^{\frac{1}{n}}}\Bigg{)}k^{\frac{1}{n}}
.\end{aligned}\end{eqnarray} Very recently, S.Y.Yolcu \cite{Y} has
improved the estimate (\ref{1.10}) to the following:
\begin{eqnarray}\label{1.11}\begin{aligned}\frac{1}{k}\sum^{k}_{j=1}\Lambda^{(\alpha)}_{j}\geq\frac{n\widetilde{C}_{n}}{n+1}Vol(D)^{-\frac{1}{n}}k^{\frac{1}{n}}
+\widetilde{M}_{n}\frac{Vol(D)^{1+\frac{1}{n}}}{Ine(D)}k^{-1/n},
\end{aligned}\end{eqnarray}where $\widetilde{C}_{n}=\frac{2\pi}{(\omega_{n})^{\frac{1}{n}}}$ and the constant $\widetilde{M}_{n}$ depends only on
the dimension $n$. Moreover, for any $\alpha\in(0,2]$, S.Y.Yolcu and
T.Yolcu \cite{YY} generalized (\ref{1.11}) as follows:
\begin{eqnarray}\label{1.12}\begin{aligned}\frac{1}{k}\sum^{k}_{j=1}\Lambda^{(\alpha)}_{j}\geq\frac{n}{n+\alpha}
\Bigg{(}\frac{2\pi}{(\omega_{n}Vol(D))^{\frac{1}{n}}}\Bigg{)}^{\alpha}k^{\frac{\alpha}{n}}
.\end{aligned}\end{eqnarray} Furthermore, S.Y.Yolcu and T.Yolcu
\cite{YY} refined the Berezin-Li-Yau inequality in the case of
fractional Laplacian $(-\Delta)^{\alpha}|_{D}$ restricted to $D$:
\begin{eqnarray}\label{1.13}\begin{aligned}\frac{1}{k}\sum^{k}_{j}\Lambda^{(\alpha)}_{j}&\geq\frac{n}{n+\alpha}
\Bigg{(}\frac{2\pi}{(\omega_{n}Vol(D))^{\frac{1}{n}}}\Bigg{)}^{\alpha}k^{\frac{\alpha}{n}}
\\&~\quad+\frac{\ell}{4(n+\alpha)}\frac{(2\pi)^{\alpha-2}}{(\omega_{n}Vol(D))^{\frac{\alpha-2}{n}}}
\frac{Vol(D)}{Ine(D)}k^{\frac{\alpha-2}{n}},
\end{aligned}\end{eqnarray}
where $\ell$ is given by $$\ell=\min\Bigg{\{}\frac{\alpha}{12},
\frac{4\alpha n\pi^{2}}{(2n+2-\alpha)\omega^{\frac{4}{n}}_{n}}
\Bigg{\}}.$$
\begin{rem}
In fact, by a direct calculation, one can check the following
inequality: $$\frac{\alpha}{12}\leq\frac{4\alpha
n\pi^{2}}{(2n+2-\alpha)\omega^{\frac{4}{n}}_{n}},$$ which
implies\begin{equation*}\begin{aligned}\frac{\ell}{4(n+\alpha)}\frac{(2\pi)^{\alpha-2}}{(\omega_{n}Vol(D))^{\frac{\alpha-2}{n}}}
\frac{Vol(D)}{Ine(D)}k^{\frac{\alpha-2}{n}}
=\frac{\alpha}{48(n+\alpha)}\frac{(2\pi)^{\alpha-2}}{(\omega_{n}Vol(D))^{\frac{\alpha-2}{n}}}
\frac{Vol(D)}{Ine(D)}k^{\frac{\alpha-2}{n}}.\end{aligned}\end{equation*}\end{rem}

The another main purpose of this paper is to provide a refinement of
the Berezin-Li-Yau type estimate. In other word, we have proved the
following:
\begin{thm}\label{thm1.2} Let $D$ be a bounded domain in an n-dimensional
Euclidean space $\mathbb{R}^{n}$. Assume that
$\Lambda^{(\alpha)}_{i}, i = 1,2,\cdots,$ is the $i$-th eigenvalue
of the fractional Laplacian $(-\Delta)^{\alpha/2}|_{D}$. Then, the
sum of its eigenvalues satisfies
\begin{eqnarray}\label{1.14}\begin{aligned}\frac{1}{k}\sum_{j=1}^{k}\Lambda^{(\alpha)}_{j}&\geq
\frac{n}{n+ \alpha}\frac{(2\pi)^{\alpha}} {(\omega_{n}Vol(
D))^{\frac{\alpha}{n}}}k^{\frac{\alpha}{n}}
\\&~\quad+\frac{\alpha}{48(n+\alpha)}\frac{(2\pi)^{\alpha-2}}{(\omega_{n}Vol(D))^{\frac{\alpha-2}{n}}}
\frac{Vol(D)}{Ine(D)}k^{\frac{\alpha-2}{n}}
\\&~\quad+\frac{\alpha(n+\alpha-2)^{2}}{\mathcal {C}(n)n(n+\alpha)^{2}}
\frac{(2\pi)^{\alpha-4}}{(\omega_{n}Vol(D))^{\frac{\alpha-4}{n}}}
\Bigg{(}\frac{Vol(D)}{Ine(D)}\Bigg{)}^{2}k^{\frac{\alpha-4}{n}}
,\end{aligned}\end{eqnarray}
where \begin{equation*}\mathcal {C}(n)=\begin{cases}4608,~~~~~~~~~~~~~~~~~~~~~~{\rm when}\ \  n\geq4,\\
6144,~~~~~~~~~{\rm when}\ \ n=2\ \ {\rm or}\ \
n=3.\end{cases}\end{equation*}In particular, the sum of its
eigenvalues satisfies
\begin{eqnarray}\label{1.15}\begin{aligned}\frac{1}{k}\sum_{j=1}^{k}\Lambda^{(2)}_{j}&\geq
\frac{nk^{\frac{2}{n}}}{n+2}\omega^{-\frac{2}{n}}_{n}(2\pi)^{2}Vol(D)^{-\frac{2}{n}}+\frac{1}{24(n+2)}\frac{Vol(D)}{Ine(D)}
\\&\quad~+\frac{nk^{-\frac{2}{n}}}{2304(n+2)^{2}}\omega_{n}^{\frac{2}{n}}(2\pi)^{-2}\Bigg{(}\frac{Vol(D)}{Ine(D)}\Bigg{)}^{2}Vol(D)^{\frac{2}{n}},
\end{aligned}\end{eqnarray}
when $\alpha=2$.
\end{thm}

\begin{rem}Observing  Theorem \ref{thm1.2}, it is not difficult to see that the coefficients (with respect to
$k^{\frac{\alpha-2}{n}}$) of the second terms in (\ref{1.14}) are
equal to that of (\ref{1.13}). In other word, we can claim that the
inequalities (\ref{1.14}) are sharper than (\ref{1.13}) since the
coefficients (with respect to $k^{\frac{\alpha-4}{n}}$) of  the
third terms in (\ref{1.14}) are positive.
\end{rem}By using Theorem \ref{thm1.2}, we can give an analogue of the Berezin-Li-Yau type
inequality for the eigenvalues of the Klein-Gordon operators
$\mathscr{H}_{0,D}$ restricted to the bounded domain $D$:
\begin{corr}Let $D$ be a bounded domain in an n-dimensional
Euclidean space $\mathbb{R}^{n}$. Assume that $\Lambda_{i}, i =
1,2,\cdots,$ is the $i$-th eigenvalue of the  Klein-Gordon operators
$\mathscr{H}_{0,D}$. Then, the sum of its eigenvalues satisfies
\begin{eqnarray}\begin{aligned}\label{1.16}\frac{1}{k}\sum_{j=1}^{k}\Lambda_{j}&\geq
\frac{n}{n+ 1}\frac{2\pi} {(\omega_{n}Vol(
D))^{\frac{1}{n}}}k^{\frac{1}{n}}
\\&~\quad+\frac{1}{48(n+1)}\frac{(2\pi)^{-1}}{(\omega_{n}Vol(D))^{-\frac{1}{n}}}
\frac{Vol(D)}{Ine(D)}k^{-\frac{1}{n}}
\\&~\quad+\frac{(n-1)^{2}}{\mathcal {C}(n)n(n+1)^{2}}
\frac{(2\pi)^{-3}}{(\omega_{n}Vol(D))^{-\frac{3}{n}}}
\Bigg{(}\frac{Vol(D)}{Ine(D)}\Bigg{)}^{2}k^{-\frac{3}{n}}
,\end{aligned}\end{eqnarray}
where \begin{equation*}\mathcal {C}(n)=\begin{cases}4608,~~~~~~~~~~~~~~~~~~~~~~{\rm when}\ \  n\geq4,\\
6144,~~~~~~~~~{\rm when}\ \ n=2\ \ {\rm or}\ \
n=3.\end{cases}\end{equation*}
\end{corr}

\section{A Key Lemma}
In order to prove the following Lemma  \ref{lem2.3} , we need the
following lemmas given by S.Y.Yolcu and T.Yolcu in \cite{YY}:

\begin{lem}\label{lem2.1}Suppose that $\varsigma: [0,\infty)\rightarrow[0,1]$ such that
$$0\leq \varsigma(s)\leq1~{\rm and}~\int^{\infty}_{0}\varsigma(s)ds=1.$$Then, there
exists $\epsilon\geq0$ such that
\begin{equation*}\int^{\epsilon+1}_{\epsilon}s^{d}ds=\int^{\infty}_{0}s^{d}\varsigma(s)ds.
\end{equation*}Moreover, we have
\begin{equation*}\int^{\epsilon+1}_{\epsilon}s^{d+\alpha}ds\leq\int^{\infty}_{0}s^{d+\alpha}\varsigma(s)ds.
\end{equation*}

\end{lem}
\begin{lem}\label{lem2.2}For $s>0,~\tau>0,~2\leq b\in\mathbb{N},~0<\alpha\leq2$, we
have the following inequality:
\begin{equation*}s^{b+\alpha}\geq\frac{b+\alpha}{b}s^{b}\tau^{\alpha}-\frac{\alpha}{b}\tau^{b+\alpha}+\frac{\alpha}{b}\tau^{b+\alpha-2}(s-\tau)^{2}.
\end{equation*}
\end{lem}

In light of Lemma \ref{lem2.1} and Lemma \ref{lem2.2}, we obtain the
following result which will play important roles in the proof of
Theorem \ref{thm1.1} and Theorem \ref{thm1.2}.

\begin{lem}\label{lem2.3}
Let $b(\geq2)$ be a positive real number and $\mu(>0)$ be defined by
(\ref{2.13}). If $\psi: [0,~+\infty)\rightarrow[0,~+\infty)$ is a
decreasing function such that $$-\mu\leq\psi^{\prime}(s)\leq0$$ and
$$A:=\int_{0}^{\infty}s^{b-1}\psi(s)ds>0,$$ then, we have
\begin{eqnarray}\label{2.1}\begin{aligned}\int^{\infty}_{0}s^{b+\alpha-1}\psi(s)ds&\geq
\frac{1}{b+\alpha}(bA)^{\frac{b+\alpha}{b}}\psi(0)^{-\frac{\alpha}{b}}\\&~\quad+\frac{\alpha}{12b(b+\alpha)\mu^{2}}(bA)^{\frac{b+\alpha-2}{b}}\psi(0)^{\frac{2b-\alpha+2}{b}}
\\&~\quad+\frac{\alpha(b+\alpha-2)^{2}}{288b^{2}(b+\alpha)^{2}\mu^{4}}(bA)^{\frac{b+\alpha-4}{b}}\psi(0)^{\frac{4b-\alpha+4}{b}},\end{aligned}\end{eqnarray}
when $b\geq4$; we have
\begin{eqnarray}\label{2.2}\begin{aligned}\int^{\infty}_{0}s^{b+\alpha-1}\psi(s)ds&\geq
\frac{1}{b+\alpha}(bA)^{\frac{b+\alpha}{b}}\psi(0)^{-\frac{\alpha}{b}}\\&~\quad+\frac{\alpha}{12b(b+\alpha)\mu^{2}}(bA)^{\frac{b+\alpha-2)}{b}}\psi(0)^{\frac{2b-\alpha+2}{b}}
\\&~\quad+\frac{\alpha(b+\alpha-2)^{2}}{384b^{2}(b+\alpha)^{2}\mu^{4}}(bA)^{\frac{b+\alpha-4}{b}}\psi(0)^{\frac{4b-\alpha+4}{b}},\end{aligned}\end{eqnarray}
when $2\leq b<4$. In particular, the inequality (\ref{2.1}) holds
when $\alpha=2$ and $b\geq2$.\end{lem}

\noindent\emph{Proof.} If we consider the following
function\begin{eqnarray*}\varrho(t)=\frac{\psi\big{(}\frac{\psi(0)}{\mu}t\big{)}}{\psi(0)},\end{eqnarray*}
then it is not difficult to see that $\varrho(0)=1$ and
$-1\leq\varrho'(t)\leq0.$ Without loss of generality, we can assume
$$\psi(0)=1~\textnormal{and}~\mu=1.$$ Define
$$E_{\alpha}:=\int^{\infty}_{0}s^{b+\alpha-1}\psi(s)ds.$$
One can assume that $E_{\alpha}<\infty$, otherwise there is nothing
to prove. By the assumption,
 we can conclude that
$$\lim_{s\rightarrow\infty}s^{b+\alpha-1}\psi(s)=0.$$ Putting
$h(s)=-\psi'(s)$ for any $s\geq0$, we get $$0\leq
h(s)\leq1~~\textnormal{and}~~\int^{\infty}_{0}h(s)ds=\psi(0)=1.$$ By
making use of integration by parts, one can get
$$\int^{\infty}_{0}s^{b}h(s)ds=b\int^{\infty}_{0}s^{b-1}\psi(s)ds=bA,$$and
$$\int^{\infty}_{0}s^{b+\alpha}h(s)ds\leq(b+\alpha)E_{\alpha},$$
since $\psi(s)>0.$ By Lemma \ref{lem2.1}, one can infer that there
exists an $\epsilon\geq0$ such that
\begin{eqnarray}\label{2.3}\int^{\epsilon+1}_{\epsilon}s^{b}ds=\int^{\infty}_{0}s^{b}h(s)ds=bA,\end{eqnarray}
and
\begin{eqnarray}\label{2.4}\int^{\epsilon+1}_{\epsilon}s^{b+\alpha}ds\leq\int^{\infty}_{0}s^{b+\alpha}h(s)ds\leq(b+\alpha)E_{\alpha}.\end{eqnarray}
Let
$$\Theta(s)=bs^{b+\alpha}-(b+\alpha)\tau^{\alpha}s^{b}+\alpha\tau^{b+\alpha}-\alpha\tau^{b+\alpha-2}(s-\tau)^{2},$$
then, by Lemma \ref{lem2.2}, we have $\Theta(s)\geq0.$ Integrating
the function $\Theta(s)$ from $\epsilon$ to $\epsilon+1$, we deduce
from (\ref{2.3}) and (\ref{2.4}), for any $\tau>0,$
\begin{eqnarray}\label{2.5}b(b+\alpha)E_{\alpha}-(b+\alpha)\tau^{\alpha}bA+\alpha\tau^{b+\alpha}\geq\frac{\alpha}{12}\tau^{b+\alpha-2}.\end{eqnarray}
Define
\begin{eqnarray}\label{2.6}f(\tau):=(b+\alpha)\tau^{\alpha}bA-\alpha\tau^{b+\alpha}+\frac{\alpha}{12}\tau^{b+\alpha-2},\end{eqnarray}
then we can obtain from (\ref{2.5}) that, for any $\tau>0,$
$$E_{\alpha}=\int^{\infty}_{0}s^{b+\alpha-1}\psi(s)ds\geq\frac{f(\tau)}{b(b+\alpha)}.$$
Taking
$$\tau=(bA)^{\frac{1}{b}}\Bigg{(}1+\frac{b+\alpha-2}{12(b+\alpha)}(bA)^{-\frac{2}{b}}\Bigg{)}^{\frac{1}{b}},$$
and substituting it into (\ref{2.6}), we obtain
\begin{eqnarray}\label{2.7}\begin{aligned}f(\tau)&=(bA)^{\frac{b+\alpha}{b}}
\Bigg{(}b-\frac{\alpha(b+\alpha-2)}{12(b+\alpha)}(bA)^{-\frac{2}{b}}\Bigg{)}\Bigg{(}1+\frac{b+\alpha-2}{12(b+\alpha)}(bA)^{-\frac{2}{b}}
\Bigg{)}^{\frac{\alpha}{b}}\\&~\quad+
\frac{\alpha}{12}(bA)^{\frac{b+\alpha-2}{b}}\Bigg{(}1+\frac{b+\alpha-2}{12(b+\alpha)}(bA)^{-\frac{2}{b}}\Bigg{)}^{\frac{b+\alpha-2}{b}}.\end{aligned}\end{eqnarray}
By using the Taylor formula, one has for $t>0$
\begin{eqnarray*}\begin{aligned}(1+t)^{\frac{\alpha}{b}}&\geq1+\frac{\alpha}{b}t+\frac{\alpha(\alpha-b)}{2b^{2}}t^{2}
+\frac{\alpha(\alpha-b)(\alpha-2b)}{6b^{3}}t^{3}\\&\quad~+\frac{\alpha(\alpha-b)(\alpha-2b)(\alpha-3b)}{24b^{4}}t^{4},
\end{aligned}\end{eqnarray*}
and
\begin{eqnarray*}\begin{aligned}(1+t)^{\frac{b+\alpha-2}{b}}&\geq
1+\frac{b+\alpha-2}{b}t+\frac{(b+\alpha-2)(\alpha-2)}{2b^{2}}t^{2}
\\&\quad~+\frac{(b+\alpha-2)(\alpha-2)(\alpha-2-b)}{6b^{3}}t^{3}
\\&\quad~+\frac{(b+\alpha-2)(\alpha-2)(\alpha-2-b)(\alpha-2-2b)}{24b^{4}}t^{4}.\end{aligned}\end{eqnarray*}
Putting
$$t=\frac{b+\alpha-2}{12(b+\alpha)}(bA)^{-\frac{2}{b}}>0,$$ one has $b-\alpha
t>0$, $\tau=(bA)^{\frac{1}{b}}(1+t)^{\frac{1}{b}},$

\begin{eqnarray}\label{2.8}\begin{aligned}&\quad~\Bigg{(}b-\frac{\alpha(b+\alpha-2)}{12(b+\alpha)}(bA)^{-\frac{2}{b}}\Bigg{)}
\Bigg{(}1+\frac{b+\alpha-2}{12(b+\alpha)}(bA)^{-\frac{2}{b}}\Bigg{)}^{\frac{\alpha}{b}}
\\&=
(b-\alpha t)(1+t)^{\frac{\alpha}{b}}\\&\geq(b-\alpha
t)\Bigg{[}1+\frac{\alpha}{b}t+\frac{\alpha(\alpha-b)}{2b^{2}}t^{2}
+\frac{\alpha(\alpha-b)(\alpha-2b)}{6b^{3}}t^{3}\\&\quad~+\frac{\alpha(\alpha-b)(\alpha-2b)(\alpha-3b)}{24b^{4}}t^{4}\Bigg{]}
\\&=
b-\frac{\alpha(\alpha+b)}{2b}t^{2}-\frac{\alpha(\alpha-b)(\alpha+b)}{3b^{2}}
t^{3}-\frac{\alpha(\alpha-b)(\alpha-2b)(\alpha+b)}{8b^{3}} t^{4}
\\&\quad~-\frac{\alpha^{2}(\alpha-b)(\alpha-2b)(\alpha-3b)}{24b^{4}}t^{5}
\\&=b-\frac{\alpha(\alpha+b)}{2b}\Bigg{(}\frac{b+\alpha-2}{12(b+\alpha)}(bA)^{-\frac{2}{b}}\Bigg{)}^{2}
\\&\quad\quad~-\frac{\alpha(\alpha-b)(\alpha+b)}{3b^{2}}
\Bigg{(}\frac{b+\alpha-2}{12(b+\alpha)}(bA)^{-\frac{2}{b}}\Bigg{)}^{3}\\&\quad\quad~-\frac{\alpha(\alpha-b)(\alpha-2b)(\alpha+b)}{8b^{3}}
\Bigg{(}\frac{b+\alpha-2}{12(b+\alpha)}(bA)^{-\frac{2}{b}}\Bigg{)}^{4}
\\&\quad\quad~-\frac{\alpha^{2}(\alpha-b)(\alpha-2b)(\alpha-3b)}{24b^{4}}\Bigg{(}\frac{b+\alpha-2}{12(b+\alpha)}(bA)^{-\frac{2}{b}}\Bigg{)}^{5}
,
\end{aligned}\end{eqnarray}
and
\begin{eqnarray}\label{2.9}\begin{aligned}&\quad~\Bigg{(}1+\frac{b+\alpha-2}{12(b+\alpha)}(bA)^{-\frac{2}{b}}\Bigg{)}^{\frac{b+\alpha-2}{b}}
\\&=(1+t)^{\frac{b+\alpha-2}{b}}\\&\geq
1+\frac{b+\alpha-2}{b}\Bigg{(}\frac{b+\alpha-2}{12(b+\alpha)}(bA)^{-\frac{2}{b}}\Bigg{)}
\\&\quad~+\frac{(b+\alpha-2)(\alpha-2)}{2b^{2}}\Bigg{(}\frac{b+\alpha-2}{12(b+\alpha)}(bA)^{-\frac{2}{b}}\Bigg{)}^{2}
\\&\quad~+\frac{(b+\alpha-2)(\alpha-2)(\alpha-2-b)}{6b^{3}}\Bigg{(}\frac{b+\alpha-2}{12(b+\alpha)}(bA)^{-\frac{2}{b}}\Bigg{)}^{3}
\\&\quad~+\frac{(b+\alpha-2)(\alpha-2)(\alpha-2-b)(\alpha-2-2b)}{24b^{4}}\Bigg{(}\frac{b+\alpha-2}{12(b+\alpha)}(bA)^{-\frac{2}{b}}\Bigg{)}^{4}
.\end{aligned}
\end{eqnarray}
Therefore, we obtain from (\ref{2.8}) and (\ref{2.9})
\begin{eqnarray}\label{2.10}\begin{aligned}f(\tau)&=(b+\alpha)\tau^{\alpha}bA-\alpha\tau^{b+\alpha}+\frac{\alpha}{12}\tau^{b+\alpha-2}
\\&\geq
(bA)^{\frac{b+\alpha}{b}}\Bigg{[}b-\frac{\alpha(\alpha+b)}{2b}\Bigg{(}\frac{b+\alpha-2}{12(b+\alpha)}(bA)^{-\frac{2}{b}}\Bigg{)}^{2}
\\&\quad~-\frac{\alpha(\alpha-b)(\alpha+b)}{3b^{2}}
\Bigg{(}\frac{b+\alpha-2}{12(b+\alpha)}(bA)^{-\frac{2}{b}}\Bigg{)}^{3}\\&\quad~-\frac{\alpha(\alpha-b)(\alpha-2b)(\alpha+b)}{8b^{3}}
\Bigg{(}\frac{b+\alpha-2}{12(b+\alpha)}(bA)^{-\frac{2}{b}}\Bigg{)}^{4}
\\&\quad~-\frac{\alpha^{2}(\alpha-b)(\alpha-2b)(\alpha-3b)}{24b^{4}}\Bigg{(}\frac{b+\alpha-2}{12(b+\alpha)}(bA)^{-\frac{2}{b}}\Bigg{)}^{5}\Bigg{]}
\\&\quad~+
\frac{\alpha}{12}(bA)^{\frac{b+\alpha-2}{b}}\Bigg{[}1+\frac{b+\alpha-2}{b}\Bigg{(}\frac{b+\alpha-2}{12(b+\alpha)}(bA)^{-\frac{2}{b}}\Bigg{)}
\\&\quad~+\frac{(b+\alpha-2)(\alpha-2)}{2b^{2}}\Bigg{(}\frac{b+\alpha-2}{12(b+\alpha)}(bA)^{-\frac{2}{b}}\Bigg{)}^{2}
\\&\quad~+\frac{(b+\alpha-2)(\alpha-2)(\alpha-2-b)}{6b^{3}}\Bigg{(}\frac{b+\alpha-2}{12(b+\alpha)}(bA)^{-\frac{2}{b}}\Bigg{)}^{3}
\\&\quad~+\frac{(b+\alpha-2)(\alpha-2)(\alpha-2-b)(\alpha-2-2b)}{24b^{4}}\Bigg{(}\frac{b+\alpha-2}{12(b+\alpha)}(bA)^{-\frac{2}{b}}\Bigg{)}^{4}\Bigg{]}
\\&=b(bA)^{\frac{b+\alpha}{b}}+\frac{\alpha}{12}(bA)^{\frac{b+\alpha-2}{b}}
+\mathcal {I}_{1}+\mathcal {I}_{2}+\mathcal {I}_{3},
\end{aligned}\end{eqnarray}where
\begin{eqnarray}\label{2.11}\mathcal{I}_{1}=\frac{\alpha(b+\alpha-2)^{2}}{288b(b+\alpha)}(bA)^{\frac{b+\alpha-4}{b}},\end{eqnarray}
\begin{eqnarray}\label{2.12}\begin{aligned}\mathcal {I}_{2}&=
\frac{\alpha(b+\alpha-2)(\alpha+2b-6)}{72b^{2}}
\Bigg{(}\frac{b+\alpha-2}{12(b+\alpha)}\Bigg{)}^{2}(bA)^{\frac{b+\alpha-6}{b}}
\\&\quad~+\frac{\alpha(b+\alpha-2)[\alpha^{2}+(5b-16)\alpha+(-6b^{2}+8b+16)]}{288b^{3}}
\\&\quad~\times\Bigg{(}\frac{b+\alpha-2}{12(b+\alpha)}\Bigg{)}^{3}(bA)^{\frac{b+\alpha-8}{b}},\end{aligned}\end{eqnarray}
\begin{eqnarray*}\mathcal {I}_{3}=\frac{\alpha\gamma}{24b^{4}}
\Bigg{(}\frac{b+\alpha-2}{12(b+\alpha)}\Bigg{)}^{5}(bA)^{\frac{b+\alpha-10}{b}},\end{eqnarray*}
and\begin{eqnarray*}\begin{aligned}\gamma&=(\alpha-2)(\alpha-2-b)(\alpha-2-2b)(b+\alpha)-\alpha(\alpha-b)(\alpha-2b)(\alpha-3b)
.\end{aligned}\end{eqnarray*} Noticing that
$$-\alpha(\alpha-b)(\alpha-2b)(\alpha-3b)\geq0,$$we have
\begin{eqnarray*}\begin{aligned}\gamma&\geq(\alpha-2)(\alpha-2-b)(\alpha-2-2b)(b+\alpha)
.\end{aligned}\end{eqnarray*}
Define$$\beta:=(\alpha-2)(\alpha-2-b)(\alpha-2-2b)(b+\alpha),$$ then
we have $\beta\leq0$ and $\gamma\geq\beta$. Therefore, we have
\begin{eqnarray}\label{2.13}\mathcal {I}_{3}\geq\frac{\alpha\beta}{24b^{4}}
\Bigg{(}\frac{b+\alpha-2}{12(b+\alpha)}\Bigg{)}^{5}(bA)^{\frac{b+\alpha-10}{b}},\end{eqnarray}

Next, we consider two cases:

\textbf{Case 1: $\textbf{b}\geq\textbf{4}$.} When $b\geq4$, for any
$\alpha\in(0,2]$, we can infer
\begin{eqnarray}\label{2.14}\begin{aligned}\alpha^{2}+(5b-16)\alpha+(-6b^{2}+8b+16)&\leq(5b-16)\alpha+(-6b^{2}+8b+20)
\\&=-6b^{2}+(8+5\alpha)b+20-16\alpha\\&\leq-24b+(8+10)b+20\\&\leq0.\end{aligned}\end{eqnarray}
Since
$(bA)^{\frac{2}{b}}\geq\frac{1}{(b+1)^{\frac{2}{b}}}\geq\frac{1}{3}$
(see \cite{CW2}), one can deduce from (\ref{2.12}) and (\ref{2.14})
\begin{eqnarray}\label{2.15}\begin{aligned}\mathcal {I}_{2}&
\geq \frac{\alpha(b+\alpha-2)(\alpha+2b-6)}{72b^{2}}
\Bigg{(}\frac{b+\alpha-2}{12(b+\alpha)}\Bigg{)}^{2}(bA)^{\frac{b+\alpha-6}{b}}
\\&\quad~+\frac{\alpha(b+\alpha-2)[\alpha^{2}+(5b-16)\alpha+(-6b^{2}+8b+16)]}{1152b^{3}}
\\&\quad~\times\Bigg{(}\frac{b+\alpha-2}{12(b+\alpha)}\Bigg{)}^{2}(bA)^{\frac{b+\alpha-6}{b}}
\\&=\frac{\alpha(b+\alpha-2)\big{[}16b(\alpha+2b-6)
+\alpha^{2}+(5b-16)\alpha+(-6b^{2}+8b+16)\big{]}}{1152b^{3}}
\\&\quad~\times\Bigg{(}\frac{b+\alpha-2}{12(b+\alpha)}\Bigg{)}^{2}(bA)^{\frac{b+\alpha-6}{b}}
\\&=\frac{\alpha(b+\alpha-2)\big{[}26b^{2}+(-88+21\alpha)b+(\alpha^{2}-16\alpha+16)\big{]}}{1152b^{3}}
\\&\quad~\times\Bigg{(}\frac{b+\alpha-2}{12(b+\alpha)}\Bigg{)}^{2}(bA)^{\frac{b+\alpha-6}{b}}.
\end{aligned}\end{eqnarray}
On the other hand,  we have
$$\mathcal {I}_{3}\geq\frac{\alpha\beta}{4608b^{4}}
\Bigg{(}\frac{b+\alpha-2}{12(b+\alpha)}\Bigg{)}^{2}(bA)^{\frac{b+\alpha-6}{b}},$$since
$\beta\leq0$ and
$(bA)^{\frac{2}{b}}\geq\frac{1}{(b+1)^{\frac{2}{b}}}\geq\frac{1}{3}$.
Therefore, the estimate of the lower bound of $\mathcal
{I}_{2}+\mathcal{I}_{3}$ can be given by
\begin{eqnarray*}\begin{aligned}\mathcal {I}_{2}+\mathcal{I}_{3}&\geq\Bigg{\{}\frac{\alpha(b+\alpha-2)
\big{[}26b^{2}+(-88+21\alpha)b+(\alpha^{2}-16\alpha+16)\big{]}}{1152b^{3}}
+\frac{\alpha\beta}{4608b^{4}}\Bigg{\}}\\&\quad~\times\Bigg{(}\frac{b+\alpha-2}{12(b+\alpha)}\Bigg{)}^{2}(bA)^{\frac{b+\alpha-6}{b}}
\\&=\frac{\alpha\{4b(b+\alpha-2)[26b^{2}+(-88+21\alpha)b+(\alpha^{2}-16\alpha+16)]+\beta\}}{4608b^{4}}
\\&\quad~\times\Bigg{(}\frac{b+\alpha-2}{12(b+\alpha)}\Bigg{)}^{2}(bA)^{\frac{b+\alpha-6}{b}}.\end{aligned}\end{eqnarray*}
Next, we will verify the following inequality
\begin{eqnarray}\label{2.16}4b(b+\alpha-2)[26b^{2}+(-88+21\alpha)b+(\alpha^{2}-16\alpha+16)]+\beta\geq0.\end{eqnarray}Indeed, since
$0<\alpha\leq2$ and $b\geq4$, we have
\begin{eqnarray}\label{2.17}\begin{aligned}&\quad~4b(b+\alpha-2)[26b^{2}+(-88+21\alpha)b+(\alpha^{2}-16\alpha+16)]+\beta
\\&=4b(b+\alpha-2)[26b^{2}+(-88+21\alpha)b+(\alpha^{2}-16\alpha+16)]
\\&\quad~+(\alpha-2)(\alpha-2-b)(\alpha-2-2b)(b+\alpha)
\\&\geq8b[26b^{2}+(-88+21\alpha)b+(\alpha^{2}-16\alpha+16)]
\\&\quad~-|(\alpha-2)(\alpha-2-b)(\alpha-2-2b)(b+\alpha)|
\\&\geq8b[26b^{2}+(-88+21\alpha)b+(\alpha^{2}-16\alpha+16)]
-2|(b+2)(2b+2)(b+2)|
\\&\geq8b[26b^{2}-88b+(\alpha^{2}-8\alpha)]-2(b+2)(2b+2)(b+2)\\&\geq8b(26b^{2}-92b) -2(b+2)(2b+2)(b+2)
\\&=204b^{3}-756b^{2}-32b-16\\&\geq60b^{2}-32b-16\\&\geq28b-16\\&\geq0.
\end{aligned}\end{eqnarray}
Thus, it is not difficult to see that the inequality (\ref{2.16})
follows from (\ref{2.17}), which implies
$$\mathcal {I}_{2}+\mathcal {I}_{3}\geq0.$$ Therefore, when
$b\geq4,$ we have
$$f(\tau)\geq b(bA)^{\frac{b+\alpha}{b}}+\frac{\alpha}{12}(bA)^{\frac{b+\alpha-2}{b}}
+\frac{\alpha(b+\alpha-2)^{2}}{288b(b+\alpha)}(bA)^{\frac{b+\alpha-4}{b}}.$$

\textbf{Case 2: $2\leq\textbf{b}<\textbf{4}$.} Uniting the equations
(\ref{2.11}), (\ref{2.12}) and (\ref{2.13}), we obtain the following
equation
\begin{eqnarray*}\begin{aligned}\mathcal {I}_{1}+\mathcal {I}_{2}+\mathcal {I}_{3}&\geq\frac{\alpha(b+\alpha-2)^{2}}{288b(b+\alpha)}(bA)^{\frac{b+\alpha-4}{b}}
\\&~\quad+\frac{\alpha(b+\alpha-2)(\alpha+2b-6)}{72b^{2}}
\Bigg{(}\frac{b+\alpha-2}{12(b+\alpha)}\Bigg{)}^{2}(bA)^{\frac{b+\alpha-6}{b}}
\\&\quad~+\frac{\alpha(b+\alpha-2)[\alpha^{2}+(5b-16)\alpha+(-6b^{2}+8b+16)]}{288b^{3}}
\\&\quad~\times\Bigg{(}\frac{b+\alpha-2}{12(b+\alpha)}\Bigg{)}^{3}(bA)^{\frac{b+\alpha-8}{b}}
\\&\quad~+\frac{\alpha\beta}{24b^{4}}
\Bigg{(}\frac{b+\alpha-2}{12(b+\alpha)}\Bigg{)}^{5}(bA)^{\frac{b+\alpha-10}{b}}
\end{aligned}\end{eqnarray*}

\begin{eqnarray*}\begin{aligned}
~\quad\quad\quad\quad&=\frac{\alpha(b+\alpha-2)^{2}}{384b(b+\alpha)}(bA)^{\frac{b+\alpha-4}{b}}
+\frac{\alpha(b+\alpha-2)^{2}}{1152b(b+\alpha)}(bA)^{\frac{b+\alpha-4}{b}}
\\&~\quad+\frac{\alpha(b+\alpha-2)\nu_{1}}{72b^{2}}
\Bigg{(}\frac{b+\alpha-2}{12(b+\alpha)}\Bigg{)}^{2}(bA)^{\frac{b+\alpha-6}{b}}
\\&\quad~+\frac{\alpha(b+\alpha-2)\nu_{2}}{288b^{3}}
\Bigg{(}\frac{b+\alpha-2}{12(b+\alpha)}\Bigg{)}^{3}(bA)^{\frac{b+\alpha-8}{b}}
\\&\quad~+\frac{\alpha\beta}{24b^{4}}
\Bigg{(}\frac{b+\alpha-2}{12(b+\alpha)}\Bigg{)}^{5}(bA)^{\frac{b+\alpha-10}{b}},\end{aligned}\end{eqnarray*}
where $$\nu_{1}:=(\alpha+2b-6),$$ and
$$\nu_{2}:=\alpha^{2}+(5b-16)\alpha+(-6b^{2}+8b+16). $$Suppose $\nu_{1}\leq0$ and $\nu_{2}\leq0$,  then we
have

\begin{eqnarray}\label{2.18}\begin{aligned}\mathcal {I}_{1}+\mathcal {I}_{2}+\mathcal {I}_{3}&\geq\frac{\alpha(b+\alpha-2)^{2}}{384b(b+\alpha)}(bA)^{\frac{b+\alpha-4}{b}}
+\frac{\alpha(b+\alpha-2)}{96b}
\Bigg{(}\frac{b+\alpha-2}{12(b+\alpha)}\Bigg{)}(bA)^{\frac{b+\alpha-4}{b}}
\\&~\quad+\frac{\alpha(b+\alpha-2)(\alpha+2b-6)}{288b^{2}}
\Bigg{(}\frac{b+\alpha-2}{12(b+\alpha)}\Bigg{)}(bA)^{\frac{b+\alpha-4}{b}}
\\&\quad~+\frac{\alpha(b+\alpha-2)[\alpha^{2}+(5b-16)\alpha+(-6b^{2}+8b+16)]}{4608b^{3}}
\\&\quad~\times\Bigg{(}\frac{b+\alpha-2}{12(b+\alpha)}\Bigg{)}(bA)^{\frac{b+\alpha-4}{b}}
\\&\quad~+\frac{\alpha\beta}{24b^{4}}
\Bigg{(}\frac{b+\alpha-2}{12(b+\alpha)}\Bigg{)}^{5}(bA)^{\frac{b+\alpha-10}{b}}
\\&=\frac{\alpha(b+\alpha-2)^{2}}{384b(b+\alpha)}(bA)^{\frac{b+\alpha-4}{b}}+\frac{\alpha(b+\alpha-2)}{96b}\mathcal {I}_{4}
\Bigg{(}\frac{b+\alpha-2}{12(b+\alpha)}\Bigg{)}(bA)^{\frac{b+\alpha-4}{b}}
\\&\quad~+\frac{\alpha\beta}{24b^{4}}
\Bigg{(}\frac{b+\alpha-2}{12(b+\alpha)}\Bigg{)}^{5}(bA)^{\frac{b+\alpha-10}{b}},
\end{aligned}\end{eqnarray}where
\begin{eqnarray*}\begin{aligned}\mathcal {I}_{4}&=1 +\frac{\alpha+2b-6}{3b}
+\frac{\alpha^{2}+(5b-16)\alpha+(-6b^{2}+8b+16)}{48b^{2}}.\end{aligned}\end{eqnarray*}Noticing
that $0<\alpha\leq2$ and $2\leq b<4$, we have
\begin{eqnarray}\label{2.19}\begin{aligned}\mathcal {I}_{4}&=\frac{48b^{2}+16b(\alpha+2b-6)+\alpha^{2}+(5b-16)\alpha+(-6b^{2}+8b+16)}{48b^{2}}
\\&=\frac{74b^{2}+(21\alpha-88)b+(\alpha^{2}-16\alpha+16)}{48b^{2}}
\\&\geq\frac{60b+21\alpha b+(\alpha^{2}-8\alpha)}{48b^{2}}\\&=\frac{60b+(21b-8)\alpha+\alpha^{2}}{48b^{2}}
\\&\geq\frac{60b}{48b^{2}}=\frac{5}{4b}.\end{aligned}\end{eqnarray}
Therefore, we derive from (\ref{2.18}) and (\ref{2.19})
\begin{eqnarray*}\begin{aligned}\mathcal {I}_{1}+\mathcal
{I}_{2}+\mathcal
{I}_{3}&\geq\frac{\alpha(b+\alpha-2)^{2}}{384b(b+\alpha)}(bA)^{\frac{b+\alpha-4}{b}}+\frac{5\alpha(b+\alpha-2)}{384b^{2}}
\Bigg{(}\frac{b+\alpha-2}{12(b+\alpha)}\Bigg{)}(bA)^{\frac{b+\alpha-4}{b}}
\\&\quad~+\frac{\alpha\beta}{24b^{4}}
\Bigg{(}\frac{b+\alpha-2}{12(b+\alpha)}\Bigg{)}^{5}(bA)^{\frac{b+\alpha-10}{b}}
\\&\geq\frac{\alpha(b+\alpha-2)^{2}}{384b(b+\alpha)}(bA)^{\frac{b+\alpha-4}{b}}
+\frac{5\alpha(b+\alpha-2)}{384b^{2}}
\Bigg{(}\frac{b+\alpha-2}{12(b+\alpha)}\Bigg{)}(bA)^{\frac{b+\alpha-4}{b}}
\\&\quad~+\frac{\alpha(b+\alpha-2)\beta}{18432(b+\alpha)b^{4}}
\Bigg{(}\frac{b+\alpha-2}{12(b+\alpha)}\Bigg{)}(bA)^{\frac{b+\alpha-4}{b}}
\\&\geq\frac{\alpha(b+\alpha-2)^{2}}{384b(b+\alpha)}(bA)^{\frac{b+\alpha-4}{b}}
\\&\quad~+\frac{\alpha(b+\alpha-2)[240b^{2}(b+\alpha)+\beta] }{18432(b+\alpha)b^{4}}
\Bigg{(}\frac{b+\alpha-2}{12(b+\alpha)}\Bigg{)}(bA)^{\frac{b+\alpha-4}{b}},\end{aligned}\end{eqnarray*}since
$(bA)^{\frac{2}{b}}\geq\frac{1}{(b+1)^{\frac{2}{b}}}\geq\frac{1}{3}$.
We define a function  $\mathcal {K}(b)$ by letting
\begin{eqnarray*}\begin{aligned}\mathcal
{K}(b)&:=240b^{2}(b+\alpha)+\beta
\\&=240b^{2}(b+\alpha)+(\alpha-2)(\alpha-2-b)(\alpha-2-2b)(b+\alpha)
,\end{aligned}\end{eqnarray*}where $b\in[2,4)$. After a direct
calculation, we have
\begin{eqnarray*}\begin{aligned}\mathcal
{K}(b)&\geq240b^{2}(b+\alpha)-|(\alpha-2)(\alpha-2-b)(\alpha-2-2b)(b+\alpha)|
\\&\geq240b^{2}(b+\alpha)-|2(2+b)(2+2b)(b+2)|
\\&\geq240b^{2}(b+\alpha)-2(2b)(3b)(2b)
\\&\geq216b^{3}+240\alpha b^{2}>0,\end{aligned}\end{eqnarray*}
which implies
\begin{eqnarray*}\begin{aligned}\mathcal {I}_{1}+\mathcal
{I}_{2}+\mathcal
{I}_{3}&\geq\frac{\alpha(b+\alpha-2)^{2}}{384b(b+\alpha)}(bA)^{\frac{b+\alpha-4}{b}}.\end{aligned}\end{eqnarray*}
For the other cases (i.e., $\nu_{1}\leq0~{\rm and}~ \nu_{2}>0$;
$\nu_{1}>0~{\rm and}~ \nu_{2}\leq0$; or $\nu_{1}>0~{\rm and}~
\nu_{2}>0$), we can also derive by using the same method
that$$\mathcal {I}_{1}+\mathcal {I}_{2}+\mathcal
{I}_{3}\geq\frac{\alpha(b+\alpha-2)^{2}}{384b(b+\alpha)}(bA)^{\frac{b+\alpha-4}{b}}.$$
Therefore, when $2\leq b\leq4,$ we have
$$f(\tau)\geq b(bA)^{\frac{b+\alpha}{b}}+\frac{\alpha}{12}(bA)^{\frac{b+\alpha-2}{b}}
+\frac{\alpha(b+\alpha-2)^{2}}{384b(b+\alpha)}(bA)^{\frac{b+\alpha-4}{b}}.$$
In particular, we can consider the case that $\alpha=2$. Noticing
that $\beta=0$ when $\alpha=2$ and $b\geq2$, we can claim that
$\mathcal {I}_{3}\geq0$. Therefore,  when $\alpha=2$  and $b\geq2$,
one can deduce
\begin{eqnarray*}\begin{aligned}\mathcal {I}_{2}+\mathcal {I}_{3}&\geq
\frac{\alpha(b+\alpha-2)(\alpha+2b-6)}{72b^{2}}
\Bigg{(}\frac{b+\alpha-2}{12(b+\alpha)}\Bigg{)}^{2}(bA)^{\frac{b+\alpha-6}{b}}
\\&\quad~+\frac{\alpha(b+\alpha-2)[\alpha^{2}+(5b-16)\alpha+(-6b^{2}+8b+16)]}{288b^{3}}
\\&\quad~\times\Bigg{(}\frac{b+\alpha-2}{12(b+\alpha)}\Bigg{)}^{3}(bA)^{\frac{b+\alpha-8}{b}}
\\&=\frac{b(b-2)}{18b^{2}}
\Bigg{(}\frac{b}{12(b+2)}\Bigg{)}^{2}(bA)^{\frac{b-4}{b}}
+\frac{-b^{2}+3b-2}{24b^{2}}
\Bigg{(}\frac{b}{12(b+2)}\Bigg{)}^{3}(bA)^{\frac{b-6}{b}}
\\&\geq\frac{b(b-2)}{18b^{2}}
\Bigg{(}\frac{b}{12(b+2)}\Bigg{)}^{2}(bA)^{\frac{b-4}{b}}
+\frac{-b^{2}+3b-2}{96b^{2}}
\Bigg{(}\frac{b}{12(b+2)}\Bigg{)}^{2}(bA)^{\frac{b-4}{b}}
\\&=\frac{13b^{2}-23b-6}{288b^{2}}
\Bigg{(}\frac{b}{12(b+2)}\Bigg{)}^{2}(bA)^{\frac{b-4}{b}}
\\&\geq\frac{26b-23b-6}{288b^{2}}
\Bigg{(}\frac{b}{12(b+2)}\Bigg{)}^{2}(bA)^{\frac{b-4}{b}}
\\&\geq0,\end{aligned}\end{eqnarray*}which implies
$$f(\tau)\geq b(bA)^{\frac{b+\alpha}{b}}+\frac{\alpha}{12}(bA)^{\frac{b+\alpha-2}{b}}
+\frac{\alpha(b+\alpha-2)^{2}}{288b(b+\alpha)}(bA)^{\frac{b+\alpha-4}{b}}.$$

This completes the proof of the Lemma \ref{lem2.3}.$$\eqno{\Box}$$

\section{Proofs of Theorem \ref{thm1.1} and Theorem \ref{thm1.2}}
In this section, we will prove the Theorem \ref{thm1.1} and Theorem
\ref{thm1.2} by using the key lemma given in section 2 (i.e., Lemma
\ref{lem2.3}).

We suppose that $D\subset\mathbb{R}^{n}$ is a bounded domain in
$\mathbb{R}^{n}$, and then its \emph{symmetric rearrangement}
$D^{\ast}$ is the open ball with the same volume as $D$,
$$D^{\ast}=\Bigg{\{}x\in\mathbb{R}^{n}\Bigg{|}~|x|<\Bigg{(}\frac{Vol(D)}{\omega_{n}}\Bigg{)}^{\frac{1}n{}}\Bigg{\}}.$$
By using a symmetric rearrangement of $D$, one can obtain
\begin{equation}\label{3.1}Ine(D)=\int_{D}|x|^{2}dx\geq\int_{D^{\ast}}|x|^{2}dx=\frac{n}{n+2}Vol(D)\Bigg{(}\frac{Vol(D)}{\omega_{n}}\Bigg{)}^{\frac{2}n{}}.\end{equation}
For the case of fractional Laplace operator, let $u_{j}^{(\alpha)}$
be an orthonormal eigenfunction corresponding to the eigenvalue
$\Lambda^{(\alpha)}_{j}$. Namely, $u_{j}^{(\alpha)}$ satisfies
\begin{equation*}\begin{cases}(-\Delta)^{\alpha/2}u^{(\alpha)}_{j} = \Lambda^{(\alpha)} u_{j}^{(\alpha)},~~~~~~~~~~~~~~~~~~~\textnormal{in}~D,\\
\int_{D}u_{i}^{(\alpha)}(x)u_{j}^{(\alpha)}(x)dx=\delta_{ij},~~~~~~~~~\textnormal{for}~
\textnormal{any}~i,j,
\end{cases}\end{equation*} where $0<\alpha\leq2$.
On the other hand, for the case of  Laplace operator, we let $v_{j}$
be an orthonormal eigenfunction corresponding to the eigenvalue
$\lambda_{j}$. Namely, $v_{j}$ satisfies
\begin{equation*}\begin{cases}\Delta v_{j} + \lambda_{j} v_{j}=0,~~~~~~~~~~~~~~~~~~~~~~~~\textnormal{in}~D,\\
v=0,~~~~~~~~~~~~~~~~~~~~~~~~~~~~~~~~~~\textnormal{on}~\partial D,\\
\int_{D}v_{i}(x)v_{j}(x)dx=\delta_{ij},~~~~~~~~~\textnormal{for}~
\textnormal{any}~i,j.
\end{cases}\end{equation*}

Thus, both $\{u^{(\alpha)}_{j}\}^{\infty}_{j=1}$ and
$\{v_{j}\}^{\infty}_{j=1}$ form an orthonormal basis of $L^{2}(D)$.
Define the functions $\varphi^{(\alpha)}_{j}$ and $\eta_{j}$ by
\begin{equation*}\varphi^{(\alpha)}_{j}(x)=\begin{cases}u^{(\alpha)}_{j}(x),~~~~~~~~~~~~~~x\in D,\\
0,~~~~~~~~~~~~~~~~x\in\mathbb{R}^{n}\backslash D,
\end{cases}\end{equation*}
and
\begin{equation*}\eta_{j}(x)=\begin{cases}v_{j}(x),~~~~~~~~~~~~~~x\in D,\\
0,~~~~~~~~~~~~~~x\in\mathbb{R}^{n}\backslash D,
\end{cases}\end{equation*}
respectively. Denote by $\widehat{\eta_{j}}(\xi)$ and
$\widehat{\varphi^{(\alpha)}_{j}}(\xi)$ the Fourier transforms of
$\eta_{j}(\xi)$ and  $\varphi^{(\alpha)}_{j}(\xi)$, then, for any
$\xi\in\mathbb{R}^{n},$ we have

$$\widehat{\varphi^{(\alpha)}_{j}}(\xi)=(2\pi)^{-n/2}\int_{\mathbb{R}^{n}}\varphi^{(\alpha)}_{j}(x)e^{i\langle x,\xi\rangle}dx
=(2\pi)^{-n/2}\int_{D}u^{(\alpha)}_{j}(x)e^{i\langle
x,\xi\rangle}dx,$$and
$$\widehat{\eta_{j}}(\xi)=(2\pi)^{-n/2}\int_{\mathbb{R}^{n}}\eta_{j}(x)e^{i\langle x,\xi\rangle}dx
=(2\pi)^{-n/2}\int_{D}v_{j}(x)e^{i\langle x,\xi\rangle}dx.$$ From
the Plancherel formula, we have
$$\int_{\mathbb{R}^{n}}\widehat{\varphi^{(\alpha)}_{i}}(x)\widehat{\varphi^{(\alpha)}_{j}}(x)dx=
\int_{\mathbb{R}^{n}}\widehat{\eta_{i}}(x)\widehat{\eta_{j}}(x)dx=\delta_{ij},$$
for any $i, j$. Since $\{u^{(\alpha)}_{j}\}^{\infty}_{j=1}$ and
$\{v_{j}\}^{\infty}_{j=1}$ are orthonormal basises in $L^{2}(D)$,
the Bessel inequality implies that
\begin{eqnarray}\label{3.2}\sum^{k}_{j=1}|\widehat{\varphi^{(\alpha)}_{j}}(\xi)|^{2}\leq(2\pi)^{-n/2}\int_{D}|e^{i\langle x,\xi\rangle}|^{2}dx
=(2\pi)^{-n/2}Vol(D),\end{eqnarray} and
\begin{eqnarray}\label{3.3}\sum^{k}_{j=1}|\widehat{\eta_{j}}(\xi)|^{2}\leq(2\pi)^{-n/2}\int_{D}|e^{i\langle
x,\xi\rangle}|^{2}dx =(2\pi)^{-n/2}Vol(D).\end{eqnarray} For
fractional Laplace operator, we observe that
\begin{eqnarray}\label{3.4}\begin{aligned}\Lambda^{(\alpha)}_{j}&=\int_{\mathbb{R}^{n}}
u^{(\alpha)}_{j}(\xi)\cdot(-\Delta)^{\alpha/2}|_{\Omega}u^{(\alpha)}_{j}(\xi)d\xi
\\&=\int_{\mathbb{R}^{n}}
u^{(\alpha)}_{j}(\xi)\cdot\mathscr{F}^{-1}[|\xi|^{\alpha}\mathscr{F}[u^{(\alpha)}_{j}(\xi)]]d\xi
\\&=\int_{\mathbb{R}^{n}}|\xi|^{\alpha}|\widehat{u^{(\alpha)}_{j}}(\xi)|^{2}d\xi,
\end{aligned}\end{eqnarray}since the support of $u^{(\alpha)}_{j}$ is
$D$ (see \cite{YY}). On the meanwhile, for the case of Laplace
operator, we have (see\cite{LiY,M})
\begin{eqnarray}\label{3.5}\lambda_{j}=\int_{\mathbb{R}^{n}}|\xi|^{2}|\widehat{v_{j}}(\xi)|^{2}d\xi.\end{eqnarray}
Since
$$\nabla\widehat{\varphi^{(\alpha)}_{j}}(\xi)=(2\pi)^{-n/2}\int_{\Omega}ixu^{(\alpha)}_{j}(x)e^{i\langle x,\xi\rangle}dx,$$
and
$$\nabla\widehat{\eta_{j}}(\xi)=(2\pi)^{-n/2}\int_{\Omega}ixv_{j}(x)e^{i\langle
x,\xi\rangle}dx,$$
we obtain
\begin{eqnarray}\label{3.6}\sum_{j=1}^{k}|\nabla\widehat{\varphi^{(\alpha)}_{j}}(\xi)|^{2}
=\sum_{j=1}^{k}|\nabla\widehat{\eta_{j}}(\xi)|^{2}=(2\pi)^{-n}\int_{\Omega}|ixe^{i\langle
x,\xi\rangle}|^{2}dx=(2\pi)^{-n}Ine(D).\end{eqnarray} Putting
$$f^{(\alpha)}(\xi):=\sum_{j=1}^{k}|\widehat{\varphi^{(\alpha)}_{j}}(\xi)|^{2},$$and
$$f(\xi):=\sum_{j=1}^{k}|\widehat{\eta_{j}}(\xi)|^{2},$$
one derives from (\ref{3.2}) and (\ref{3.3}) that $0\leq
f^{(\alpha)}(\xi)\leq(2\pi)^{-n}Vol(D)$ and $0\leq
f(\xi)\leq(2\pi)^{-n}Vol(D)$, it follows from (\ref{3.6}) and the
Cauchy-Schwarz inequality that
\begin{eqnarray*}\begin{aligned}|\nabla f^{(\alpha)}(\xi)|&\leq2\Bigg{(}\sum_{j=1}^{k}|\widehat{\varphi^{(\alpha)}_{j}}(\xi)|^{2}\Bigg{)}^{1/2}
\Bigg{(}\sum_{j=1}^{k}|\nabla\widehat{\varphi^{(\alpha)}_{j}}(\xi)|^{2}\Bigg{)}^{1/2}\\&\leq2(2\pi)^{-n}\sqrt{Ine(D)Vol(D)},
\end{aligned}\end{eqnarray*}and
\begin{eqnarray*}\begin{aligned}|\nabla f(\xi)|&\leq2\Bigg{(}\sum_{j=1}^{k}|\widehat{\eta_{j}}(\xi)|^{2}\Bigg{)}^{1/2}
\Bigg{(}\sum_{j=1}^{k}|\nabla\widehat{\eta_{j}}(\xi)|^{2}\Bigg{)}^{1/2}\\&\leq2(2\pi)^{-n}\sqrt{Ine(D)Vol(D)},
\end{aligned}\end{eqnarray*}for
every $\xi\in\mathbb{R}^{n}$. Furthermore, by using (\ref{3.4}) and
(\ref{3.5}), we have
\begin{eqnarray}\label{3.7}\begin{aligned}\sum^{k}_{j=1}\Lambda^{(\alpha)}_{j}&=
\sum^{k}_{j=1}\int_{\mathbb{R}^{n}}|\xi|^{\alpha}|\widehat{u^{(\alpha)}_{j}}(\xi)|^{2}d\xi
=\sum^{k}_{j=1}\int_{\mathbb{R}^{n}}|\xi|^{\alpha}|\widehat{\varphi^{(\alpha)}_{j}}(\xi)|^{2}d\xi
\\&=\int_{\mathbb{R}^{n}}|\xi|^{\alpha}f^{(\alpha)}(\xi)d\xi,\end{aligned}\end{eqnarray}and
\begin{eqnarray}\label{3.8}\begin{aligned}\sum^{k}_{j=1}\lambda_{j}&=
\sum^{k}_{j=1}\int_{\mathbb{R}^{n}}|\xi|^{2}|\widehat{v_{j}}(\xi)|^{2}d\xi
=\sum^{k}_{j=1}\int_{\mathbb{R}^{n}}|\xi|^{2}|\widehat{\eta_{j}}(\xi)|^{2}d\xi
\\&=\int_{\mathbb{R}^{n}}|\xi|^{2}f(\xi)d\xi.\end{aligned}\end{eqnarray}
From the Parseval's identity, we derive
\begin{eqnarray}\label{3.9}\begin{aligned}\int_{\mathbb{R}^{n}}f^{(\alpha)}(\xi)d\xi&=\sum_{j=1}^{k}\int_{\mathbb{R}^{n}}|\widehat{\varphi^{(\alpha)}_{j}}(x)|^{2}dx
=\sum_{j=1}^{k}\int_{D}|\widehat{u^{(\alpha)}_{j}}(x)|^{2}dx
\\&=\sum_{j=1}^{k}\int_{D}|u^{(\alpha)}_{j}(x)|^{2}dx=k.\end{aligned}\end{eqnarray}Similarly,
we have \cite{CW1,M}
\begin{eqnarray}\label{3.10}\begin{aligned}\int_{\mathbb{R}^{n}}f(\xi)d\xi&=\sum_{j=1}^{k}\int_{\mathbb{R}^{n}}|\widehat{\eta_{j}}(x)|^{2}dx
=\sum_{j=1}^{k}\int_{D}|\widehat{v_{j}}(x)|^{2}dx
\\&=\sum_{j=1}^{k}\int_{D}|v_{j}(x)|^{2}dx=k.\end{aligned}\end{eqnarray}
Let $h$ be a nonnegative bounded continuous function on $D$ and
$h^{\ast}$ is its symmetric decreasing rearrangement, then we have
(see \cite{Ba,CW1})
\begin{eqnarray}\label{3.11}\begin{aligned}\int_{\mathbb{R}^{n}} h(x)dx = \int_{\mathbb{R}^{n}} h^{\ast}(x)dx
= n\omega_{n
}\int_{0}^{\infty}s^{n-1}g(s)ds\end{aligned}\end{eqnarray} and
\begin{eqnarray}\label{3.12}\begin{aligned}\int_{\mathbb{R}^{n}}|x|^{\alpha}h(x)dx \geq
\int_{\mathbb{R}^{n}}|x|^{\alpha}h^{\ast}(x)dx = n\omega_{n
}\int_{0}^{\infty}s^{n+\alpha-1}g(s)ds,\end{aligned}\end{eqnarray}
where $\alpha\in (0,2]$ and $g(|x|)=h^{\ast}(x)$. Putting
 $\delta:=\sup|\nabla h|,$ then we can obtain
\begin{eqnarray}\label{3.13}\begin{aligned}-\delta\leq g'(s)\leq0\end{aligned}\end{eqnarray} for almost every $s$.
More detail information on symmetric decreasing rearrangements will
be found in \cite{Ba,CW1,PS}.

To be brief, we will drop the superscript $\alpha$ to denote
$f^{(\alpha)}$ by $f_{1}$ and let $f_{2}=f$. Assume that
$f_{i}^{\ast}$ is the symmetric decreasing rearrangement of $f_{i}$
($i=1,2$), according to (\ref{3.9}), (\ref{3.10}) and (\ref{3.11}),
we have
\begin{eqnarray}\label{3.14}k=\int_{\mathbb{R}^{n}}f_{i}(\xi)d\xi=\int_{\mathbb{R}^{n}}f_{i}^{\ast}
(\xi)d\xi=n\omega_{n}\int^{\infty}_{0}s^{n-1}\phi_{i}(s)ds,\end{eqnarray}
where  $\phi_{i}(x)=f_{i}^{\ast}(|x|)$ and $i=1,2$.

Applying the symmetric decreasing rearrangement to $f_{i}$, and
noting that
\begin{eqnarray}\label{3.15}\delta_{i}\leq2(2\pi)^{-n}\sqrt{Ine(\Omega)Vol(\Omega)}:=\sigma,\end{eqnarray}where
$\delta_{i}=\sup|\nabla f_{i}|$, we obtain from (\ref{3.13})
$$-\sigma\leq-\delta_{i}\leq
\phi_{i}'(s)\leq0,$$where $i=1,2$. By(\ref{3.1}), we have
$$\sigma\geq2(2\pi)^{-n}(\frac{n}{n+2})^{\frac{1}{2}}\omega_{n}^{-\frac{1}{n}}
Vol(D)^{\frac{n+1}{n}}
\geq(2\pi)^{-n}\omega_{n}^{-\frac{1}{n}}Vol(D)^{\frac{n+1}{n}},$$since
$n\geq2$. Moreover, by using (\ref{3.7}), (\ref{3.8}) and
(\ref{3.12}), we have
\begin{eqnarray}\label{3.16}\begin{aligned}\sum^{k}_{j=1}\Lambda^{(\alpha)}_{j}&=
\int_{\mathbb{R}^{n}}|\xi|^{\alpha}f^{(\alpha)}(\xi)d\xi=
\int_{\mathbb{R}^{n}}|\xi|^{\alpha}f_{1}(\xi)d\xi
\\&\geq\int_{\mathbb{R}^{n}}|\xi|^{\alpha}f^{\ast}_{1}(\xi)d\xi
=n\omega_{n}\int^{\infty}_{0}s^{n+\alpha-1}\phi_{1}(\xi)d\xi,\end{aligned}\end{eqnarray}and
\begin{eqnarray}\label{3.17}\begin{aligned}\sum^{k}_{j=1}\lambda_{j}&=
\int_{\mathbb{R}^{n}}|\xi|^{2}f(\xi)d\xi=
\int_{\mathbb{R}^{n}}|\xi|^{2}f_{2}(\xi)d\xi
\\&\geq\int_{\mathbb{R}^{n}}|\xi|^{2}f^{\ast}_{2}(\xi)d\xi
=n\omega_{n}\int^{\infty}_{0}s^{n+1}\phi_{2}(\xi)d\xi.\end{aligned}\end{eqnarray}

\emph{Proof of Theorem 1.1.} In order to apply Lemma \ref{lem2.3},
from (\ref{3.14}), (\ref{3.15}) and the definition of $A$, we take
$$b=n,\ \ \psi(s)=\phi_{2}(s),\ \ A=\frac{k}{n\omega_{n}},\ \ {\rm and}\
\ \mu=\sigma=2(2\pi)^{-n}\sqrt{Vol(D)Ine(D)}.$$ Therefore, we can
obtain from Lemma \ref{lem2.3} and (\ref{3.17}) that
\begin{eqnarray*}\begin{aligned}\sum^{k}_{j=1}\lambda_{j}&=n\omega_{n}E_{2}\\&\geq
\frac{n\omega_{n}}{n(n+2)}\Bigg{[}n(nA)^{\frac{n+2}{n}}\phi_{2}(0)^{-\frac{2}{n}}+\frac{nA\phi_{2}(0)^{2}}
{6\sigma^{2}}+\frac{n(nA)^{\frac{n-2}{n}}}{144(n+2)\sigma^{4}}\phi_{2}(0)^{\frac{4n+2}{n}}\Bigg{]}
\\&=
\frac{\omega_{n}}{n+2}\Bigg{[}n\Bigg{(}\frac{k}{\omega_{n}}\Bigg{)}^{\frac{n+2}{n}}t^{-\frac{2}{n}}+\frac{\frac{k}{\omega_{n}}t^{2}}
{6\sigma^{2}}+\frac{n(\frac{k}{\omega_{n}})^{\frac{n-2}{n}}}{144(n+2)\sigma^{4}}t^{\frac{4n+2}{n}}\Bigg{]}
\\&=\frac{n}{n+2}\omega_{n}^{-\frac{2}{n}}k^{\frac{n+2}{n}}t^{-\frac{2}{n}}+\frac{kt^{2}}
{6(n+2)\sigma^{2}}+\frac{n\omega_{n}^{\frac{2}{n}}k^{\frac{n-2}{n}}}{144(n+2)^{2}\sigma^{4}}t^{\frac{4n+2}{n}},
\end{aligned}\end{eqnarray*}
where $t=\phi_{2}(0)$. Let
$$F(t)=\frac{n}{n+2}\omega_{n}^{-\frac{2}{n}}k^{\frac{n+2}{n}}t^{-\frac{2}{n}}+\frac{kt^{2}}
{6(n+2)\sigma^{2}}+\frac{n\omega_{n}^{\frac{2}{n}}k^{\frac{n-2}{n}}}{144(n+2)^{2}\sigma^{4}}t^{\frac{4n+2}{n}},$$
then one can has
$$F'(t)=-\frac{2}{n+2}\omega_{n}^{-\frac{2}{n}}k^{\frac{n+2}{n}}t^{-\frac{n+2}{n}}+\frac{kt}
{3(n+2)\sigma^{2}}+\frac{4n+2}{144(n+2)^{2}\sigma^{4}}\omega_{n}^{\frac{2}{n}}k^{\frac{n-2}{n}}t^{\frac{3n+2}{n}}.$$
Since $F'(t)$ is increasing on $(0,(2\pi)^{-n}Vol(D)]$, then it is
easy to see that $F(t)$ is decreasing on $(0,(2\pi)^{-n}Vol(D)]$ if
$F'((2\pi)^{-n}Vol(D))<0$. Indeed,
\begin{eqnarray*}\begin{aligned}F'((2\pi)^{-n}Vol(D))&\leq-\frac{2}{n+2}\omega_{n}^{-\frac{2}{n}}k^{\frac{n+2}{n}}
((2\pi)^{-n}Vol(D))^{-\frac{n+2}{n}}\\&~\quad+\frac{k((2\pi)^{-n}Vol(D))}
{3(n+2)\Big{[}(2\pi)^{-n}\omega_{n}^{-\frac{1}{n}}Vol(D)^{\frac{n+1}{n}}\Big{]}^{2}}
\\&~\quad+\frac{(4n+2)\omega_{n}^{\frac{2}{n}}k^{\frac{n-2}{n}}((2\pi)^{-n}Vol(D))^{\frac{3n+2}{n}}}
{144(n+2)^{2}\Big{[}(2\pi)^{-n}\omega_{n}^{-\frac{1}{n}}Vol(D)^{\frac{n+1}{n}}\Big{]}^{4}}
\\&=-\frac{2}{n+2}(2\pi)^{n+2}\omega_{n}^{-\frac{2}{n}}k^{\frac{n+2}{n}}Vol(D)^{-\frac{n+2}{n}}
\\&~\quad+\frac{1}{3(n+2)}(2\pi)^{n}\omega_{n}^{\frac{2}{n}}k
Vol(D)^{-\frac{n+2}{n}}\\&~\quad+\frac{4n+2}{144(n+2)^{2}}(2\pi)^{n-2}\omega_{n}^{\frac{6}{n}}k^{\frac{n-2}{n}}Vol(D)^{-\frac{n+2}{n}}
\\&=\frac{(2\pi)^{n}k}{n+2}\omega_{n}^{\frac{2}{n}}Vol(D)^{-\frac{n+2}{n}}\mathcal {J}
,
\end{aligned}\end{eqnarray*} where
\begin{eqnarray*}\begin{aligned}
\mathcal {J}&=\frac{1}{3}+\frac{4n+2}
{144(n+2)}(2\pi)^{-2}\omega_{n}^{\frac{4}{n}}k^{-\frac{2}{n}}-2(2\pi)^{2}k^{\frac{2}{n}}\omega_{n}^{-\frac{4}{n}}
\\&<\frac{1}{3}+\frac{4(n+2)}
{144(n+2)}(2\pi)^{-2}\omega_{n}^{\frac{4}{n}}-2(2\pi)^{2}\omega_{n}^{-\frac{4}{n}}
\\&=\frac{1}{3}+\frac{1}
{36}(2\pi)^{-2}\omega_{n}^{\frac{4}{n}}-2(2\pi)^{2}\omega_{n}^{-\frac{4}{n}}
\\&<\frac{1}{3}+\frac{1}{72}-4
\\&<0,\end{aligned}
\end{eqnarray*}which implies that $F'((2\pi)^{-n}Vol(D))<0.$
Here, we use the inequality
$\frac{\omega_{n}^{\frac{4}{n}}}{(2\pi)^{2}}<\frac{1}{2}$. We can
replace $\phi_{2}(0)$ by $(2\pi)^{-n}Vol(D)$ to obtain
\begin{eqnarray*}\frac{1}{k}\sum_{j=1}^{k}\lambda_{j}&\geq&
\frac{nk^{\frac{2}{n}}}{n+2}\omega^{-\frac{2}{n}}_{n}(2\pi)^{2}Vol(D)^{-\frac{2}{n}}+\frac{1}{24(n+2)}\frac{Vol(D)}{Ine(D)}
\nonumber\\&&+\frac{nk^{-\frac{2}{n}}}{2304(n+2)^{2}}\omega_{n}^{\frac{2}{n}}(2\pi)^{-2}
\Bigg{(}\frac{Vol(D)}{Ine(D)}\Bigg{)}^{2}Vol(D)^{\frac{2}{n}}.\end{eqnarray*}
since $\sigma=2(2\pi)^{-n}\sqrt{Vol(D)Ine(D)}.$

This completes the proof of Theorem \ref{thm1.1}.$$\eqno{\square}$$

Next, we will give the proof of Theorem \ref{thm1.2}.\\

\emph{Proof of Theorem \ref{thm1.2}.}: Define the function
$\phi_{1}(x)$   by $\phi_{1}(|x|):=f_{1}^{\ast}(x)$. Then we know
that $\phi_{1}:[0,~+\infty)\rightarrow[0,~(2\pi)^{-n}V(\Omega)]$ is
a non-increasing function with respect to $|x|$. Taking
$$b=n,\ \ \psi(s)=\phi_{1}(s),\ \ A=\frac{k}{n\omega_{n}},\ \ {\rm and}\
\ \mu=\sigma=2(2\pi)^{-n}\sqrt{Vol(D)Ine(D)},$$ we can obtain from
Lemma \ref{lem2.3} and (\ref{3.16}) that
\begin{eqnarray}\label{3.18}\begin{aligned}\sum_{j=1}^{k}\Lambda^{(\alpha)}_{j}&\geq n\omega_{n}\int^{\infty}_{0}s^{n+\alpha-1}\phi_{1}(s)ds\\&\geq
\frac{n\omega_{n}\Big{(}\frac{k}{\omega_{n}}\Big{)}^{\frac{n+\alpha}{n}}}{n+\alpha}\phi_{1}(0)^{-\frac{\alpha}{n}}
+\frac{\alpha
\omega_{n}\Big{(}\frac{k}{\omega_{n}}\Big{)}^{\frac{n+\alpha-2}{n}}}{12(n+\alpha)\sigma^{2}}\phi_{1}(0)^{\frac{2n-\alpha+2}{n}}
\\&~\quad+\frac{\alpha(n+\alpha-2)^{2}\omega_{n}\Big{(}\frac{k}{\omega_{n}}\Big{)}^{\frac{n+\alpha-4}{n}}}{\mathcal {C}_{1}(n)n(n+\alpha)^{2}\sigma^{4}}
\phi_{1}(0)^{\frac{4n-\alpha+4}{n}},\end{aligned}\end{eqnarray}
where \begin{equation*}\mathcal {C}_{1}(n)=\begin{cases}288,~~~~~~~~~~~~~~~~~~~~~~{\rm when}\ \  n\geq4,\\
384,~~~~~~~~{\rm when}\ \ n=2\ \ {\rm or}\ \
n=3.\end{cases}\end{equation*} Moreover, we define a function
$\xi(t)$ by letting
\begin{eqnarray}\label{3.19}\begin{aligned}\xi(t)&=
\frac{n\omega_{n}}{n+\alpha}\Bigg{(}\frac{k}{\omega_{n}}\Bigg{)}^{\frac{n+\alpha}{n}}t^{-\frac{\alpha}{n}}
+\frac{\alpha\omega_{n}}{12(n+\alpha)\sigma^{2}}\Bigg{(}\frac{k}{\omega_{n}}\Bigg{)}^{\frac{n+\alpha-2}{n}}t^{\frac{2n-\alpha+2}{n}}
\\&~\quad+\frac{\alpha(n+\alpha-2)^{2}\omega_{n}}{\mathcal {C}_{1}(n)n(n+\alpha)^{2}\sigma^{4}}
\Bigg{(}\frac{k}{\omega_{n}}\Bigg{)}^{\frac{n+\alpha-4}{n}}t^{\frac{4n-\alpha+4}{n}}.\end{aligned}\end{eqnarray}
Differentiating (\ref{3.19}) with respect to the variable $t$, it is
not difficult to see that

\begin{eqnarray}\label{3.20}\begin{aligned}\xi^{\prime}(t)&=
\frac{\alpha\omega_{n}}{n+\alpha}\Bigg{(}\frac{k}{\omega_{n}}\Bigg{)}^{\frac{n+\alpha}{n}}t^{-\frac{\alpha}{n}-1}\Bigg{[}-1
+\frac{(2n-\alpha+2)}{12n\sigma^{2}}\Bigg{(}\frac{k}{\omega_{n}}\Bigg{)}^{-\frac{2}{n}}t^{\frac{2n+2}{n}}
\\&~\quad+\frac{(4n-\alpha+4)(n+\alpha-2)^{2}}{\mathcal {C}_{1}(n)n^{2}(n+\alpha)\sigma^{4}}
\Bigg{(}\frac{k}{\omega_{n}}\Bigg{)}^{-\frac{4}{n}}t^{\frac{4n+4}{n}}\Bigg{]}.
\end{aligned}\end{eqnarray} Letting \begin{eqnarray}\label{3.21}\begin{aligned}
\zeta(t)=\xi^{\prime}(t)(\frac{n+\alpha}{\alpha\omega_{n}})(\frac{k}{\omega_{n}})^{-\frac{n+\alpha}{n}}t^{\frac{\alpha}{n}+1},\end{aligned}\end{eqnarray}
and noting that
$\sigma\geq(2\pi)^{-n}\omega_{n}^{-\frac{1}{n}}Vol(D)^{\frac{n+1}{n}}$,
we can obtain from (\ref{3.20}) and (\ref{3.21}) that
\begin{eqnarray}\label{3.22}\begin{aligned}
\zeta(t)&= -1
+\frac{(2n-\alpha+2)}{12n\sigma^{2}}\Bigg{(}\frac{k}{\omega_{n}}\Bigg{)}^{-\frac{2}{n}}t^{\frac{2n+2}{n}}
\\&~\quad+\frac{(4n-\alpha+4)(n+\alpha-2)^{2}}{\mathcal {C}_{1}(n)n^{2}(n+\alpha)\sigma^{4}}
\Bigg{(}\frac{k}{\omega_{n}}\Bigg{)}^{-\frac{4}{n}}t^{\frac{4n+4}{n}}
\\&\leq-1
+\frac{(2n-\alpha+2)}{12n(2\pi)^{-2n}\omega_{n}^{-\frac{2}{n}}Vol(D)^{\frac{2(n+1)}{n}}}
\Bigg{(}\frac{k}{\omega_{n}}\Bigg{)}^{-\frac{2}{n}}t^{\frac{2n+2}{n}}
\\&~\quad+\frac{(4n-\alpha+4)(n+\alpha-2)^{2}}{\mathcal {C}_{1}(n)n^{2}(n+\alpha)(2\pi)^{-4n}\omega_{n}^{-\frac{4}{n}}Vol(D)^{\frac{4(n+1)}{n}}}
\Bigg{(}\frac{k}{\omega_{n}}\Bigg{)}^{-\frac{4}{n}}t^{\frac{4n+4}{n}}.\end{aligned}\end{eqnarray}
It is easy to see that the right hand side of (\ref{3.22}) is an
increasing function of $t$. Therefore, if the right hand side of
(\ref{3.22}) is less than $0$ when we take $t=(2\pi)^{-n}Vol(D)$,
which is equivalent to say that
\begin{eqnarray}\label{3.23}\begin{aligned}
\zeta(t)&\leq-1
+\frac{(2n-\alpha+2)}{12n}k^{-\frac{2}{n}}\frac{\omega_{n}^{\frac{4}{n}}}{(2\pi)^{2}}
\\&~\quad+\frac{(4n-\alpha+4)(n+\alpha-2)^{2}}{\mathcal {C}_{1}(n)n^{2}(n+\alpha)}
k^{-\frac{4}{n}}\frac{\omega_{n}^{\frac{8}{n}}}{(2\pi)^{4}}
\\&\leq0,\end{aligned}\end{eqnarray}we can claim from
(\ref{3.23}) that $\xi'(t)\leq0$ on $(0,(2\pi)^{-n}V(\Omega)].$ By a
direct calculation, we can obtain
\begin{eqnarray}\label{3.24}\begin{aligned}
\zeta(t)&\leq-1 +\frac{(2n-\alpha+2)}{12n}
+\frac{(4n-\alpha+4)(n+\alpha-2)^{2}}{\mathcal
{C}_{1}(n)n^{2}(n+\alpha)}\\&\leq-1 +\frac{(2n+n)}{12n}
+\frac{(4n+2n)(n+n)^{2}}{\mathcal {C}_{1}(n)n^{3}}
\\&=-\frac{3}{4}+\frac{24}{\mathcal {C}_{1}(n)}
\\&\leq0,\end{aligned}\end{eqnarray}
since $\frac{\omega_{n}^{\frac{4}{n}}}{(2\pi)^{2}}<1.$ Thus, it is
easy to see from (\ref{3.21}) and (\ref{3.24}) that $\xi'(t)\leq0$,
which implies that $\xi(t)$ is a decreasing function on
$(0,(2\pi)^{-n}Vol(D)].$

On the other hand, we notice that
$0<\phi_{1}(0)\leq(2\pi)^{-n}Vol(D)$ and right hand side of the
formula (\ref{3.18}) is $\xi(\phi_{1}(0))$, which is a decreasing
function of $\phi_{1}(0)$ on $(0,(2\pi)^{-n}Vol(D)]$. Therefore,
$\phi_{1}(0)$ can be replaced by $(2\pi)^{-n}Vol(D)$ in (\ref{2.1})
which gives the following inequality:
\begin{eqnarray*}\begin{aligned}\frac{1}{k}\sum_{j=1}^{k}\Lambda^{(\alpha)}_{j}&\geq
\frac{n}{n+ \alpha}\frac{(2\pi)^{\alpha}} {(\omega_{n}Vol(
D))^{\frac{\alpha}{n}}}k^{\frac{\alpha}{n}}
\\&~\quad+\frac{\alpha}{48(n+\alpha)}\frac{(2\pi)^{\alpha-2}}{(\omega_{n}Vol(D))^{\frac{\alpha-2}{n}}}
\frac{Vol(D)}{Ine(D)}k^{\frac{\alpha-2}{n}}
\\&~\quad+\frac{\alpha(n+\alpha-2)^{2}}{\mathcal {C}(n)n(n+\alpha)^{2}}
\frac{(2\pi)^{\alpha-4}}{(\omega_{n}Vol(D))^{\frac{\alpha-4}{n}}}
\Bigg{(}\frac{Vol(D)}{Ine(D)}\Bigg{)}^{2}k^{\frac{\alpha-4}{n}}
,\end{aligned}\end{eqnarray*}
where \begin{equation*}\mathcal {C}(n)=\begin{cases}4608,~~~~~~~~~~~~~~~~~~~~~~{\rm when}\ \  n\geq4,\\
6144,~~~~~~~~~{\rm when}\ \ n=2\ \ {\rm or}\ \
n=3.\end{cases}\end{equation*}

In particular, when $\alpha=2$, we can get the inequality
(\ref{1.15}) by using the same method as the proof of Theorem
\ref{thm1.1}.

This completes the proof of Theorem \ref{thm1.2}.
$$\eqno{\Box}$$
\begin{ack}\emph{The authors wish to express their gratitude to Prof. Q.-M. Cheng for
continuous encouragement and enthusiastic help.}\end{ack}

\underline{}

\begin{flushleft}
\medskip\noindent

Guoxin Wei, School of Mathematical Sciences, South China Normal
University, 510631, Guangzhou, China,
weigx03@mails.tsinghua.edu.cn\end{flushleft}

\begin{flushleft}
\medskip\noindent

He-Jun Sun,  Department of Applied Mathematics, College of Science,
Nanjing University of Science and Technology, 210094,  Nanjing,
China hejunsun@163.com
\end{flushleft}

\begin{flushleft}
\medskip\noindent

Lingzhong Zeng, Department of Mathematics, Graduate School of
Science and Engineering, Saga University, Saga 840-8502, Japan,
lingzhongzeng@yeah.net

\end{flushleft}
\end{document}